\documentclass[11pt]{article} % for LaTex2e

\usepackage{amsfonts}
\usepackage{latexsym}
\catcode`\@=11\relax
\newwrite\@unused
\def\typeout#1{{\let\protect\string\immediate\write\@unused{#1}}}

\def\psglobal#1{
\immediate\special{ps: plotfile #1 }}
\def\psfiginit{\typeout{psfiginit}
\immediate\psglobal{figtex.pro}%
\special{ps:: /TeXMagnification {\the\mag} def}%marco:look @ figtex.pro.
}

\def\@nnil{\@nil}
\def\@empty{}
\def\@psdonoop#1\@@#2#3{}
\def\@psdo#1:=#2\do#3{\edef\@psdotmp{#2}\ifx\@psdotmp\@empty \else
    \expandafter\@psdoloop#2,\@nil,\@nil\@@#1{#3}\fi}
\def\@psdoloop#1,#2,#3\@@#4#5{\def#4{#1}\ifx #4\@nnil \else
       #5\def#4{#2}\ifx #4\@nnil \else#5\@ipsdoloop #3\@@#4{#5}\fi\fi}
\def\@ipsdoloop#1,#2\@@#3#4{\def#3{#1}\ifx #3\@nnil
       \let\@nextwhile=\@psdonoop \else
      #4\relax\let\@nextwhile=\@ipsdoloop\fi\@nextwhile#2\@@#3{#4}}
\def\@tpsdo#1:=#2\do#3{\xdef\@psdotmp{#2}\ifx\@psdotmp\@empty \else
    \@tpsdoloop#2\@nil\@nil\@@#1{#3}\fi}
\def\@tpsdoloop#1#2\@@#3#4{\def#3{#1}\ifx #3\@nnil
       \let\@nextwhile=\@psdonoop \else
      #4\relax\let\@nextwhile=\@tpsdoloop\fi\@nextwhile#2\@@#3{#4}}
\def\psdraft{
	\def\@psdraft{0}
	\def\@psdraftspecial{100}
}
\def\psdraftspecial{
	\def\@psdraft{0}
	\def\@psdraftspecial{0}
}
\def\psfull{
	\def\@psdraft{100}
}
\psfull

\newif\if@prologfile
\newif\if@postlogfile
\newif\if@bbllx
\newif\if@bblly
\newif\if@bburx
\newif\if@bbury
\newif\if@height
\newif\if@width
\newif\if@rheight
\newif\if@rwidth
\newif\if@clip
\newif\if@right
\newif\if@left
\newif\if@toplines
\newif\if@box
\newif\if@caption
\newif\if@surround
\newif\if@captionwidth
\newif\if@captionwrite
\newif\if@captionopen
\def\@p@@sclip#1{\@cliptrue}
\def\@p@@sfile#1{%\typeout{file is #1}
		\def\@p@sfile{#1}
}
\def\@p@@sfigure#1{
		\def\@p@sfile{#1}
}
\def\@p@sfake{\hbox to 0pt{\hss Whatever\hss}}
\def\@p@@sbbllx#1{
		\@bbllxtrue
		\@d@mscratch=#1
		\edef\@p@sbbllx{\number\@d@mscratch}
}
\def\@p@@sbblly#1{
		\@bbllytrue
		\@d@mscratch=#1
		\edef\@p@sbblly{\number\@d@mscratch}
}
\def\@p@@sbburx#1{
		\@bburxtrue
		\@d@mscratch=#1
		\edef\@p@sbburx{\number\@d@mscratch}
}
\def\@p@@sbbury#1{
		\@bburytrue
		\@d@mscratch=#1
		\edef\@p@sbbury{\number\@d@mscratch}
}
\def\@p@@sheight#1{
		\@heighttrue
		\@d@mscratch=#1
   		\edef\@p@sheight{\number\@d@mscratch}
}
\def\@p@@swidth#1{
		\@widthtrue
		\@d@mscratch=#1
		\edef\@p@swidth{\number\@d@mscratch}
}
\def\@p@@srheight#1{
		\@rheighttrue
		\@d@mscratch=#1
		\edef\@p@srheight{\number\@d@mscratch}
}
\def\@p@@srwidth#1{
		\@rwidthtrue
		\@d@mscratch=#1
		\edef\@p@srwidth{\number\@d@mscratch}
}
\def\@p@@sright#1{\@righttrue \@surroundtrue}
\def\@p@@sleft#1{\@lefttrue \@surroundtrue}
\def\@p@@sextraheight#1{\@d@mextraheight=#1}
\def\@p@@sbox#1{\@boxtrue}
\def\@p@@scaption#1{\@captiontrue}
\def\@p@@stoplines#1{
		\@toplinestrue
		\@c@ttoplines=#1
}
\def\@p@@scaptionwidth#1{
		\@captionwidthtrue
	  	\@d@mcaptionwidth=#1
}
\def\@p@@scaptionwrite#1{
		\global\@captionwritetrue
		\global\@w@rname=\expandafter{\jobname_captions.tex}
		\typeout{Captions are written to \the\@w@rname}
}

\def\@p@@sprolog#1{\@prologfiletrue\def\@prologfileval{#1}}
\def\@p@@spostlog#1{\@postlogfiletrue\def\@postlogfileval{#1}}
\def\@cs@name#1{\csname #1\endcsname}
\def\@setparms#1=#2,{\@cs@name{@p@@s#1}{#2}}
\def\ps@init@parms{
		\@bbllxfalse \@bbllyfalse
		\@bburxfalse \@bburyfalse
		\@heightfalse \@widthfalse
		\@rheightfalse \@rwidthfalse
		\def\@p@sbbllx{}\def\@p@sbblly{}
		\def\@p@sbburx{}\def\@p@sbbury{}
		\def\@p@sheight{}\def\@p@swidth{}
		\def\@p@srheight{}\def\@p@srwidth{}
		\def\@p@sfile{}
		\def\@p@scost{10}
		\def\@sc{}
		\@prologfilefalse
		\@postlogfilefalse
		\@clipfalse
		\@rightfalse \@leftfalse
		\@boxfalse \@captionfalse
		\@toplinesfalse \@surroundfalse
		\@d@mextraheight=0pt
 		\@c@ttoplines=0
		\@pshape={} \def\@p@srheight@total{}
		\@captionwidthfalse \@d@mcaptionwidth=0pt
}
\def\parse@ps@parms#1{
	 	\@psdo\@psfiga:=#1\do
		   {\expandafter\@setparms\@psfiga,}}
\newif\ifno@bb
\newif\ifnot@eof
\newread\ps@stream
\newtoks\@linetok
\def\bb@missing{
	\typeout{psfig: searching \@p@sfile \space  for bounding box}
	\openin\ps@stream=\@p@sfile
	\no@bbtrue
	\not@eoftrue
	\catcode`\%=12
	\loop
		\read\ps@stream to \line@in
		\global\@linetok=\expandafter{\line@in}
		\ifeof\ps@stream \not@eoffalse \fi
		\@bbtest{\@linetok}
		\if@bbmatch\not@eoffalse\expandafter\bb@cull\the\@linetok\fi
	\ifnot@eof \repeat
	\catcode`\%=14
}	
\catcode`\%=12
\newif\if@bbmatch
\def\@bbtest#1{\expandafter\@a@\the#1%%BoundingBox:\@bbtest\@a@}
\long\def\@a@#1%%BoundingBox:#2#3\@a@{
     \ifx\@bbtest#2\@bbmatchfalse\else\@bbmatchtrue\fi}
\long\def\bb@cull#1 #2 #3 #4 #5 {
	\@d@mscratch=#2 bp\edef\@p@sbbllx{\number\@d@mscratch}
	\@d@mscratch=#3 bp\edef\@p@sbblly{\number\@d@mscratch}
	\@d@mscratch=#4 bp\edef\@p@sbburx{\number\@d@mscratch}
	\@d@mscratch=#5 bp\edef\@p@sbbury{\number\@d@mscratch}
	\no@bbfalse
}
\catcode`\%=14
\def\compute@bb{
		\no@bbfalse
		\if@bbllx \else \no@bbtrue \fi
		\if@bblly \else \no@bbtrue \fi
		\if@bburx \else \no@bbtrue \fi
		\if@bbury \else \no@bbtrue \fi
		\ifno@bb \bb@missing \fi
		\ifno@bb \typeout{FATAL ERROR: no bb supplied or found}
			\no-bb-error
		\fi
		\count203=\@p@sbburx
		\count204=\@p@sbbury
		\advance\count203 by -\@p@sbbllx
		\advance\count204 by -\@p@sbblly
		\edef\@bbw{\number\count203}
		\edef\@bbh{\number\count204}
}
\def\in@hundreds#1#2#3{\count240=#2 \count241=#3
		     \count100=\count240	% 100 is first digit #2/#3
		     \divide\count100 by \count241
		     \count101=\count100
		     \multiply\count101 by \count241
		     \advance\count240 by -\count101
		     \multiply\count240 by 10
		     \count101=\count240	%101 is second digit of #2/#3
		     \divide\count101 by \count241
		     \count102=\count101
		     \multiply\count102 by \count241
		     \advance\count240 by -\count102
		     \multiply\count240 by 10
		     \count102=\count240	% 102 is the third digit
		     \divide\count102 by \count241
		     \count200=#1\count205=0
		     \count201=\count200
			\multiply\count201 by \count100
		     	\advance\count205 by \count201
		     \count201=\count200
			\divide\count201 by 10
		     	\multiply\count201 by \count101
			\advance\count205 by \count201
		     \count201=\count200
			\divide\count201 by 100
			\multiply\count201 by \count102
			\advance\count205 by \count201
		     \edef\@result{\number\count205}
}
\def\compute@wfromh{
		\in@hundreds{\@p@sheight}{\@bbw}{\@bbh}
		\edef\@p@swidth{\@result}
}
\def\compute@hfromw{
		\in@hundreds{\@p@swidth}{\@bbh}{\@bbw}
		\edef\@p@sheight{\@result}
}
\def\compute@handw{
		\if@height
			\if@width
			\else
				\compute@wfromh
			\fi
		\else
			\if@width
				\compute@hfromw
			\else
				\edef\@p@sheight{\@bbh}
				\edef\@p@swidth{\@bbw}
			\fi
		\fi
}
\def\compute@resv{
		\if@rheight \else \edef\@p@srheight{\@p@sheight} \fi
		\if@rwidth \else \edef\@p@srwidth{\@p@swidth} \fi
		\edef\@p@srheight@total{\@p@srheight}
}
\newtoks\@pshape
\def\@c@ttoplines{\count120}
\def\@c@theightcount{\count121}
\def\@c@tshapecount{\count122}
\newdimen\@d@mwidthshape
\newdimen\@d@mextraheight
\newdimen\@d@mscratch
\def\compute@parshape{
	\if@right
		\if@left
	   		\typeout{error: Can't have both left and right set}
			\@leftfalse
		\fi
	\fi
	\@d@mscratch=\@p@swidth truesp
	\divide \@d@mscratch by 19
	\multiply \@d@mscratch by 20
	\edef\@p@swidthdimen{\the\@d@mscratch}
	\@c@tshapecount=\@c@ttoplines
 	\@d@mscratch=\baselineskip
	\multiply \@d@mscratch by \@c@ttoplines
	\advance \@d@mscratch by .4\baselineskip
    	\edef\@p@stopdistance{\the\@d@mscratch }
	\@d@mscratch=\@p@sheight truesp
	\divide \@d@mscratch by 2
	\edef\@p@shalfboxheight{\the\@d@mscratch}
	\if@toplines
		\loop \@pshape=\expandafter{\the\@pshape 0pt \hsize}
		\advance\@c@ttoplines by -1
		\ifnum\@c@ttoplines>0 \repeat
	\fi
   	\@c@theightcount=\@p@srheight@total
	\advance \@c@theightcount by \@d@mextraheight
	\divide  \@c@theightcount by \baselineskip
	\advance \@c@theightcount by 1
    	\advance \@c@tshapecount by \@c@theightcount
	\advance \@c@theightcount by -1
	\@d@mwidthshape=\hsize
     	\advance \@d@mwidthshape by -\@p@swidthdimen
	\if@left
		\def\@moveshape{0pt}
		\ifnum\@c@theightcount>0
		      	\loop
			\@pshape=%
			\expandafter{\the\@pshape %
					\@p@swidthdimen \@d@mwidthshape}
			\advance \@c@theightcount by -1
			\ifnum\@c@theightcount>0 \repeat
		\else
			\advance \@c@tshapecount by 1
		\fi
	\fi
	\if@right
		\@d@mscratch=\hsize
		\advance \@d@mscratch by -\@p@swidth truesp
		\edef\@moveshape{\@d@mscratch}
		\ifnum\@c@theightcount>0
			\loop
			\@pshape=\expandafter{\the\@pshape 0pt \@d@mwidthshape}
			\advance \@c@theightcount by -1
			\ifnum\@c@theightcount>0 \repeat
		\else
			\advance \@c@tshapecount by 1
		\fi
	\fi
	\ifnum \@p@srheight=0
		\@pshape={}
		\@c@tshapecount = 0
	\else
	 	\@pshape=\expandafter{\the\@pshape 0pt \hsize}
	\fi
}
\def\@p@ssetsurroundboxes{
	\global\parshape=\count122 \the\@pshape		% \@c@tshapecount
 	\moveright\@moveshape
	\vbox to 0pt\bgroup\hskip0pt\vskip\@p@stopdistance
}
\newtoks\@captiontok
\newbox\@b@xcaption
\newdimen\@d@mcaptionwidth
\newdimen\@d@mcaptionheight
\newwrite\@w@rcaption
\newtoks\@w@rname
\def\setcaption#1{\@captiontok={#1}}
\def\@set@caption{
	\setbox\@b@xcaption\vbox{\hsize\@d@mcaptionwidth
	\tolerance=9000 \baselineskip=11.4pt
	\noindent\relax\the\@captiontok}
	\if@captionwrite
		\if@captionopen
		\else
			\global\@captionopentrue
			\immediate\openout\@w@rcaption=\the\@w@rname
		\fi
		\immediate\write\@w@rcaption{\the\@captiontok}
		\immediate\write\@w@rcaption{}
	\fi
}
\def\compute@caption{
	\if@captionwidth
	\else
		\@d@mcaptionwidth = \@p@swidth truesp
		\divide \@d@mcaptionwidth by 20
		\multiply \@d@mcaptionwidth by 17
	\fi
	\@set@caption
	\@d@mcaptionheight=\ht\@b@xcaption
	\if@rheight
	\else
		\count100=\@p@srheight
	   	\advance \count100 by \@d@mcaptionheight
	   	\advance \count100 by \bigskipamount
	   	\advance \count100 by \medskipamount
	   	\edef\@p@srheight@total{\number\count100}
	\fi
}
\newif\if@alreadyjtem \@alreadyjtemfalse
\def\newpar{\hfil\vadjust{\vskip\parskip}%
	{\count100=\parskip
	\count101=\baselineskip
	\divide\count101 by 10  \multiply\count101 by 3
	\advance \count100 by \count101
	\divide\count100 by \baselineskip
	\advance\count100 by \prevgraf
	\global\prevgraf=\count100}%
	\break\if@alreadyjtem\else\indent\fi%
}
\let\sav@par=\par
\def\jtem#1{%
    	\if@alreadyjtem\else\bgroup\fi
	\def\par{\sav@par\egroup\sav@par}
	\if@alreadyjtem\else\leftskip\parindent\fi
	\@alreadyjtemtrue
	\noindent\hskip0pt
	\llap{#1\ }\ignorespaces
}
\def\compute@sizes{%
	\compute@bb
	\compute@handw
  	\compute@resv
	\if@caption
		\compute@caption
	\fi
	\if@surround
		\compute@parshape
	\fi
}
\def\@p@sdospecials{
	\ifnum\@p@scost<\@psdraft
	       	\typeout{psfig: including \@p@sfile \space }
	\fi
	\special{ps::[begin] 	\@p@swidth \space \@p@sheight \space
			\@p@sbbllx \space \@p@sbblly \space
			\@p@sbburx \space \@p@sbbury \space
			startTexFig \space }
	\ifnum\@p@scost<\@psdraft
		\if@clip
			\typeout{(clip)}
			\special{ps:: \@p@sbbllx \space \@p@sbblly \space
				\@p@sbburx \space \@p@sbbury \space
			    	doclip \space }
		\fi
	\fi
	\if@box
		\typeout{(box)}
  		\special{ps:: \@p@sbbllx \space \@p@sbblly \space
			\@p@sbburx \space \@p@sbbury \space
		    	dobox \space }
	\fi
	\ifnum\@p@scost<\@psdraft
		\if@prologfile
	    		\special{ps: plotfile \@prologfileval \space }
		\fi
		\special{ps: plotfile \@p@sfile \space }
    		\if@postlogfile
			\special{ps: plotfile \@postlogfileval \space }
		\fi
	\fi
	\special{ps::[end] endTexFig \space }
}
\newif\if@putinvbox

\def\psfig#1{{%
	\ifhmode%
		\vbox\bgroup
		\@putinvboxtrue
	\else
		\@putinvboxfalse
	\fi
       	\ps@init@parms
	\parse@ps@parms{#1}
       	\compute@sizes
	\if@surround
		\psfig@for@surround
	\else
		\psfig@for@regular
	\fi
	\if@putinvbox
       		\egroup
	\fi
}}
\def\psfig@for@surround{%
	\@p@ssetsurroundboxes
	\ifnum\@p@scost<\@psdraft
		\@p@sdospecials
		\vbox to \@p@srheight true sp{\vss}
       	\else
		\if@box
			\@p@sdospecials
		\fi
		\vbox to \@p@srheight true sp{
			\vskip\@p@shalfboxheight
			\hbox to \@p@srwidth true sp{
				\hss
				\ifnum\@p@scost<\@psdraftspecial
					\@p@sfile
				\else
					\@p@sfake
				\fi
      				\hss
			}
		\vss
		}
	\fi
	\if@caption
		\medskip
		\hbox to \@p@srwidth true sp{
			\hss
			\box\@b@xcaption
			\hss
		}
 		\medskip
	\fi
	\vss\egroup
	\vskip-\parskip
}

\def\psfig@for@regular{%
	\if@putinvbox
	\else
		\vskip\parskip
	\fi
	\ifnum\@p@scost<\@psdraft
		\@p@sdospecials
		\vbox to \@p@srheight true sp{%
			\hbox to \@p@srwidth true sp{
			\hfil
			}
		\vfil
		}
       	\else
		\if@box
			\@p@sdospecials
		\fi
	    	\vbox to \@p@srheight true sp{
			\vss
			\hbox to \@p@srwidth true sp{
				\hss
				\ifnum\@p@scost<\@psdraftspecial
					\@p@sfile
				\else
					\@p@sfake
				\fi
				\hss
			}
		    	\vss
		}
	\fi
	\if@caption
		\medskip
		\hbox to \@p@srwidth true sp{
			\hss
			\box\@b@xcaption
			\hss
		}
		\bigskip
	\fi
	\if@putinvbox
	\else
		\vskip-\parskip
	\fi
}
\catcode`\@=12\relax

\psfiginit

\font\scriptsizebbfont=msbm7 scaled \magstep 1
 1
\font\footnotesizebbfont=msbm9 scaled \magstep 0
 2
\font\bbfont=msbm9 scaled \magstep1  % 1.2 pt
 1
 2
 3
 4

\def\scriptsizeBbb#1{\hbox{\scriptsizebbfont #1}}

\def\footnotesizeBbb#1{\hbox{\footnotesizebbfont #1}}

\def\Bbb#1{\hbox{\bbfont #1}}

\newcommand{\Bl}{\mbox{\it Bl}}

\newcommand{\Eq}{\mbox{\rm Eq}}

\newcommand{\NE}{\mbox{\it NE}\,}

\newcommand{\Spec}{\mbox{\it Spec}\,}

\newcommand{\scriptsizealg}{\mbox{\it\scriptsize alg}}
\newcommand{\bboxtimes}{\Box\hspace{-1.75ex}\raisebox{.15ex}{$\times$}\,}

\newcommand{\ev}{\mbox{\it ev}\,}

\newcommand{\rel}{\mbox{\it\scriptsize rel}}

\newcommand{\virt}{\mbox{\scriptsize\it virt}}

\begin{document}

\enlargethispage{23cm}

\begin{titlepage}

$ $

\vspace{-2cm} % Re: -1.5cm for PC; -2.5cm for UT-Math-system

\noindent\hspace{-1cm}
\parbox{6cm}{\small July 2004}\
   \hspace{6.5cm}\
   \parbox{5cm}{math.AG/0408147}

\vspace{1.5cm}
%\vspace{2cm}

%title
\centerline{\large\bf
 A degeneration formula of Gromov-Witten invariants
}
\vspace{1ex}
\centerline{\large\bf
 with respect to a curve class for degenerations from blow-ups
}
%\vspace{1ex}
%\centerline{\large\bf
% ??????????
%} % end-title

\vspace{1.5cm}
%\vspace{2cm}
%authors-'n-addresses
\centerline{\large
  Chien-Hao Liu
  \hspace{1ex} and \hspace{1ex}
  Shing-Tung Yau
}

\vspace{3em}

%abstract%
\begin{quotation}
\centerline{\bf Abstract}
\vspace{0.3cm}
\baselineskip 12pt  %13pt for [12pt] style
{\small
 In two very detailed, technical, and fundamental works,
  Jun Li constructed
   a theory of Gromov-Witten invariants for a singular scheme
    of the gluing form $Y_1\cup _D Y_2$ that arises from
    a degeneration $W/{\Bbb A}^1$ and
   a theory of relative Gromov-Witten invariants for
    a codimension-$1$ relative pair $(Y,D)$.
 As a summit, he derived a degeneration formula that relates
  a finite summation of the usual Gromov-Witten invariants of
  a general smooth fiber $W_t$ of $W/{\Bbb A}^1$ to the Gromov-Witten
  invariants of the singular fiber $W_0=Y_1\cup_D Y_2$ via gluing
  the relative pairs $(Y_1,D)$ and $(Y_2,D)$.
 The finite sum mentioned above depends on a relative ample line
  bundle $H$ on $W/{\Bbb A}^1$.
 His theory has already applications to string theory and mathematics
  alike.
 For other new applications of Jun Li's theory, one needs a refined
  degeneration formula that depends on a curve class $\beta$ in
  $A_{\ast}(W_t)$ or $H_2(W_t;{\Bbb Z})$, rather than on the line
  bundle $H$.
 Some monodromy effect has to be taken care of to deal with this.
 For the simple but useful case of a degeneration $W/{\Bbb A}^1$ that
  arises from blowing up a trivial family $X\times{\Bbb A}^1$,
  we explain how the details of Jun Li's work can be employed to
  reach such a desired degeneration formula.
 The related set $\Omega_{(g,k;\beta)}$ of admissible triples adapted
  to $(g,k;\beta)$ that appears in the formula can be obtained
  via an analysis on the intersection numbers of relevant cycles and
  a study of Mori cones that appear in the problem.
 This set is intrinsically determined by $(g,k;\beta)$ and the normal
  bundle ${\cal N}_{Z/X}$ of the smooth subscheme $Z$ in $X$ to be
  blown up.
} %endsmall
\end{quotation}

%\bigskip
\vspace{7.6em}

\baselineskip 12pt
{\footnotesize
} %endfootnotesize

\medskip

\noindent {\small
MSC number 2000$\,$:
 14N35, 81T30.
} % end-small

\bigskip

\baselineskip 11pt
{\footnotesize
\noindent{\bf Acknowledgements.}
 We thank
  Jun Li
 for a conversation, communications, and comments on this work.
 His work influences the current work obviously. 
 C.-H.L.\ would like to thank in addition 
  Shiraz Minwalla, Mircea Mustata, Mihnea Popa, and Andrew Strominger
 for many valuable lectures and conversations;
  and Ling-Miao Chou for the tremendous moral support.
 The work is supported by NSF grants DMS-9803347 and DMS-0074329.
} %endfootnotesize

\end{titlepage}

%paper
\newpage
$ $

\vspace{-4em}  % Re: -4cm for PC; -6cm for UT-Math-system

%short heading
\centerline{\sc
 A Degeneration Formula of Gromov-Witten Invariants
}

\vspace{2em}

\baselineskip 14pt  %Re: 14pt for [11pt] style
                    %Re: 15pt for [12pt] style.

\begin{flushleft}
{\Large\bf 0. Introduction and outline.}
\end{flushleft}

\begin{flushleft}
{\bf Introduction.}
\end{flushleft}
The understanding of string world-sheet instantons, their moduli space,
 and the exact computation of the string correlations functions have
 been important problems in string theory,
 e.g.\ [A-D-K-M-V], [A-G-N-T1], [A-G-N-T2], [A-K-M-V], [A-M],
       [B-C-O-V], [C-dlO-G-P], [O-V], [Wi], and [Y-Y].
The A-model discussions of string theorists' work are formulated in part
 as the Gromov-Witten theory on the mathematical side.
Readers are referred to [Li1: Sec.\ 0] and [Li2: Sec.\ 0]
 for a quick review of the development and the fundamental literatures
 of the Gromov-Witten theory.
In [Li1] and [Li2], J.\ Li  constructed a Gromov-Witten theory for
 a singular variety of the gluing form $Y_1\cup_D Y_2$ that arises from
 a degeneration $W\rightarrow {\Bbb A}^1$ of smooth projective varieties.
Along the way as necessary ingredients, he constructed
 the stack ${\mathfrak W}$ of expanded degenerations associated to
  $W\rightarrow {\Bbb A}^1$ [Li1: Sec.\ 1],
 the stack ${\mathfrak M}({\mathfrak W},\Gamma)$ of stable morphisms
  of topological type $\Gamma$ from curves to the universal degeneration
  over ${\mathfrak W}$ [Li1: Sec.\ 2 and Sec.\ 3],
 the stack ${\mathfrak Y}_i^{\rel}$ of expanded relative pairs
  associated to $(Y_i,D)$ [Li1: Sec.\ 4],
 the stack ${\mathfrak M}({\mathfrak Y}_i^{\rel},\Gamma_i)$ of relative
  stable morphisms [Li1: Sec.\ 4],
 a perfect obstruction theory for ${\mathfrak M}({\mathfrak W},\Gamma)$
  and ${\mathfrak M}({\mathfrak Y}^{\rel}_i,\Gamma_i)$
  [Li2: Sec.\ 1, Sec.\ 2, Sec.\ 5]
  that gives rise to virtual fundamental classes
  $[{\mathfrak M}({\mathfrak W},\Gamma)]^{\virt}$ and
  $[{\mathfrak M}({\mathfrak Y}_i^{\rel},\Gamma_i)]^{\virt}$, and
 various distinguished Cartier divisors:
  $({\mathbf L}_0, {\mathbf t}_0)$ and
  $({\mathbf L}_{\eta},{\mathbf s}_{\eta})$'s on the stack
  ${\mathfrak M}({\mathfrak W},\Gamma)$ [Li2: Sec.\ 3.1].
These are used to define
 Gromov-Witten invariants of $W_0:=Y_1\cup_DY_2$ with values
  in ${\Bbb Q}$ and
 relative Gromov-Witten classes for relative pairs $(Y_i,D)$
  with values in $H_{\ast}(D^r;{\Bbb Q})$.
By construction, these invariants of the singular $W_0$ have the nice
 constant behavior under smoothening of $W_0$.
As a summit of his work, these invariants are linked together in
 a degeneration formula that relates (a summation of) the ordinary
 Gromov-Witten invariants of a smooth fiber $W_t$ of the degeneration
 $W\rightarrow {\Bbb A}^1$ to those of the singular fiber $W_0$
 via gluing the associated relative pairs $(Y_1,D)$ and $(Y_2,D)$
 [Li2: Sec.\ 3, Sec.\ 4]:
 $$
  \Psi_{\Gamma}^{W_t}(\alpha(t),\zeta)\;
   =\; \sum_{\eta\in\overline{\Omega}_{\Gamma}}\;
      \frac{{\mathbf m}(\eta)}{|\Eq(\eta)|}\,
        \sum_{j\in K_{\eta}}\,
         \left[ \Psi_{\Gamma_1}^{Y_1^{\rel}}(j_1^{\ast}\alpha(0),
                  \zeta_{\eta,1,j})\,
                \bullet\,
                \Psi_{\Gamma_2}^{Y_2^{\rel}}(j_2^{\ast}\alpha(0),
                  \zeta_{\eta,2,j})
         \right]_0\,,
 $$
[Li2: Sec.\ 0 and Sec.\ 3.2].
The theory he developed already has important applications to
 string theory and mathematics alike,
 e.g.\ [B-P], [L-S], and [G-V], [L-L-Z1], [L-L-Z2].

In J.\ Li's formula above, the left-hand side of the identity involves
 a summation over the curve classes represented by lattice points on
 a compact slice in the Mori-cone $\overline{\NE}(W_t)$ (determined
 by the degree of a fixed relative ample line bundle $H$ on $W/{\Bbb A}^1$)
 of the usual Gromov-Witten invariants with respect to classes in
 $A_1(W_t)$ or $H_2(W_t;{\Bbb Z})$.
For some new applications of J.\ Li's theory, one needs a similar
 degeneration formula whose both sides depend only on a given curve class
 in $A_1(W_t)$ or $H_2(W_t;{\Bbb Z})$,
 rather than on a relative ample line bundle $H$ on $W/{\Bbb A}^1$.
In the current work, for the special type of degenerations that arise
 from blowing up a trivial family $X\times{\Bbb A}^1$ over ${\Bbb A}^1$,
 we derive one such formula from J.\ Li's work:
 (Sec.\ 3: Theorem 3.3)
 $$
  \Psi_{(g,k;\beta)}^X(\alpha,\zeta)\;
   =\; \sum_{\eta\in\overline{\Omega}_{(g,k;\beta)}}\;
       \frac{{\mathbf m}(\eta)}{|\Eq(\eta)|}\,
         \sum_{j\in K_{\eta}}\,
          \left[ \Psi_{\Gamma_1}^{Y_1^{\rel}}(j_1^{\ast}\alpha(0),
                   \zeta_{\eta,1,j})\,
                 \bullet\,
                 \Psi_{\Gamma_2}^{Y_2^{\rel}}(j_2^{\ast}\alpha(0),
                   \zeta_{\eta,2,j})
          \right]_0\,,
 $$
 where
  $\Psi_{(g,k;\beta)}^X(\alpha,\zeta)$ on the left-hand side is the usual
   Gromov-Witten invariants defined via the moduli stack
   $\overline{\cal M}_{g,k}(X,\beta)$ with $\beta\in A_1(X)$ or
   $H_2(X;{\Bbb Z})$ and
  the expression on the right-hand side depends only on the normal bundle
   ${\cal N}_{Z/X}$ of the smooth subscheme $Z\subset X$ to be blown up
   and the triple $(g,k;\beta)$.

\bigskip

\noindent
{\bf Convention.}
 This work follows the notations and the terminology of [Li1] and [Li2]
  closely, except where confusions may occur.
 Other notations follow [Ha], [Fu], [De], and [K-M].
 All schemes are over ${\Bbb C}$ and all points are referred to closed
  points.

\bigskip

\begin{flushleft}
{\bf Outline.}
\end{flushleft}
{\small
\baselineskip 11pt  %13pt
\begin{quote}
 1. J.\ Li's degeneration formula of Gromov-Witten invariants.

 2. A degeneration formula with respect to a curve class.

 3. The $H$-(in)dependence of $\Omega_{(g,k;\beta)}^H$.
\end{quote}
} %endsmall

\bigskip

\baselineskip 14pt  %Re: 14pt for [11pt] style
                    %Re: 15pt for [12pt] style.

\bigskip

\section{J.\ Li's degeneration formula of GW invariants.}

Since [Li1] and [Li2] are very detailed and a summary of his work is
 already given in [G-V: Sec.\ 2] with many insights, we will recall
 here only definitions that are most relevant and needed to the current
 work.

Let $\pi:(W,W_0)\rightarrow ({\Bbb A}^1,{\mathbf 0})$ be a degeneration
 with smooth fiber $W_t$ for $t\ne {\mathbf 0}$ and a degenerate fiber
 $W_0$ of the gluing form $Y_1\cup_D Y_2$ over a point
 ${\mathbf 0}\in{\Bbb A}^1$,
 where $Y_i$ are smooth varieties with a smooth divisor $D_i\simeq D$.
Associated to $W/{\Bbb A}^1$ is the Artin stack ${\mathfrak W}$ of
 expanded degenerations [Li1: Sec.\ 1].
${\mathfrak W}$ with its universal family are the descent of
 the local standard models of expanded degenerations
 $W[n]\rightarrow {\Bbb A}^1[n]:= {\Bbb A}^{n+1}$ constructed from
 $\pi:W\rightarrow {\Bbb A}^1$, [Li1: Sec.\ 1.1].
A {\it stable morphism} from a prestable curve ${\cal C}/S$
 to the local model $W[n]/{\Bbb A}^1[n]$ is a diagram:
 (in notation $f:{\cal C}/S\rightarrow W[n]/{\Bbb A}^1[n]$)
 $$
  \begin{array}{ccc}
   {\cal C}    & \stackrel{f}{\longrightarrow}   & W[n]            \\
   \downarrow  &                                 & \downarrow      \\
   S           &  \rightarrow                    & {\Bbb A}^1[n]
  \end{array}
 $$
 such that $f$ is nondegenerate, predeformable, and that at every closed
 point on $S$, $f_s$ has only a finite automorphism group,
 ([Li1: Sec.\ 2.2, Sec.\ 3.1]).
Fix a relative ample line bundle $H$ on $W/{\Bbb A}^1$ and let
 $(g,k;d)$ be a triple of integers, then 
 the space of all stable morphisms from prestable curves to the universal
 family of ${\mathfrak W}$ that have arithmetic genus $g$, $k$ marked points
 on the domain and whose induced image on fibers of $W/{\Bbb A}^1$ has
 constant $H$-degree $d$ is a Deligne-Mumford stack
 ${\mathfrak M}({\mathfrak W},(g,k;d))$,
 [Li1: Theorem 3.10].

For a relative pair $(Y,D)$, where both $Y$ and $D$ are smooth and
 $D$ is a divisor in $Y$, J.\ Li constructed also a stack
 ${\mathfrak Y}^{\rel}$ of expanded relative pairs, [Li1: Sec.\ 4].
Its local model $Y[n]/{\Bbb A}^n$ has a distinguished divisor
 $D[n]/{\Bbb A}^n$ induced from $D\subset Y$, [Li1: Sec.\ 4.1].
There is a tautological morphism
 $\varphi: (Y[n],D[n])\rightarrow (Y,D)$ by construction.

\bigskip

\noindent
{\bf Definition 1.1 [admissible weighted graph].} ([Li1: Definition 4.6].)
 {\rm
  Given a relative pair $(Y,D)$, an {\it admissible weighted graph}
   $\Gamma$ for $(Y,D)$ is a graph without edges together
   with the following data:
  \begin{itemize}
   \item [(1)]
    an ordered collection of {\it legs},
    an ordered collection of weighted {\it roots}, and
    two {\it weight functions} on the vertex set
     $g:V(\Gamma)\rightarrow {\Bbb Z}_{\ge 0}$ and
     $b:V(\Gamma)\rightarrow A_1(Y)/\sim_{\scriptsizealg}$,

   \item[(2)]
    $\Gamma$ is {\it relatively connected} in the sense that either
     $|V(\Gamma)|=1$ or each vertex in $V(\Gamma)$ has at least one root
     attached to it.
  \end{itemize}
} % end-definition

\bigskip

For a fixed admissible weighted graph $\Gamma$ for $(Y,D)$, J.\ Li
 defines similarly a {\it relative stable morphism} of type $\Gamma$
 to the universal family of ${\mathfrak Y}^{\rel}$,
 [Li1: Definition 4.7. Definition 4.8].
They are locally modelled on a diagram:
 (in notation $f:{\cal C}/S\rightarrow Y[n]/{\Bbb A}^n$)
 $$
  \begin{array}{ccc}
   {\cal C}    & \stackrel{f}{\longrightarrow}   & Y[n]            \\
   \downarrow  &                                 & \downarrow      \\
   S           &  \rightarrow                    & {\Bbb A}^n
  \end{array}\,,
 $$
 in which $f$ satisfies the similar non-degenerate, predeformable, and
 stable conditions with the extra condition:
 $f^{-1}(D[n])=$ the distinguished divisor on ${\cal C}/S$
  adapted to the weighted roots of $\Gamma$.
The space of all such morphisms is a Deligne-Mumford stack
 ${\mathfrak M}({\mathfrak Y}^{\rel},\Gamma)$,
 [Li1: Definition 4.9, Theorem 4.10].
Suppose that $\Gamma$ has $k$-many legs and $r$-many roots,
 then there is an evaluation map
 $$
  \ev\,:\; {\mathfrak M}({\mathfrak Y}^{\rel},\Gamma)\;
     \longrightarrow\; Y^k
 $$
 associated to the ordinary $k$ marked points on the domain curve
 of a relative stable morphism and a distinguished evaluation map
 $$
  {\mathbf q}\;:\; {\mathfrak M}({\mathfrak Y}^{\rel},\Gamma)\;
     \longrightarrow\; D^{\,r}
 $$
 associated to the $r$ distinguished marked points that are required
 to be the only points that are mapped to $D[n]$ in a local model.

\bigskip

\noindent
{\bf Definition 1.2 [admissible triple].}
([Li1: Definition 4.11].)
{\rm
 Given a gluing $Y_1\cup_D Y_2$ of relative pairs,
 let $\Gamma_1$ and $\Gamma_2$ be a pair of admissible weighted graphs
  for $(Y_1,D)$ and $(Y_2,D)$ respectively.
 Suppose that $\Gamma_1$ and $\Gamma_2$ have identical number $r$
  of roots and $k_1$-many and $k_2$-many legs respectively.
 Let $k=k_1+k_2$ and $I\subset \{1,\,\ldots,\, k\}$ be a set of $k_1$
  elements.
 Then $(\Gamma_1,\Gamma_2,I)$ is called an {\it admissible triple} if
  the following conditions hold:
  \begin{itemize}
   \item[(1)]
    the weights on the roots of $\Gamma_1$ and $\Gamma_2$ coincide:
     $\mu_{1,i}=\mu_{2,i}$, $i=1,\,\ldots,\, r\,$;

   \item[(2)]
    after connecting the $i$-th root of $\Gamma_1$ and the $i$-th root
    of $\Gamma_2$ for all $i$, the resulting new graph with $k$ legs
    and no roots is connected.
  \end{itemize}
} % end-definition

\bigskip

Given an admissible triple $\eta=(\Gamma_1,\Gamma_2,I)$ as above
 with $Y_1\cup _D Y_2=$ the degenerate fiber $W_0$ of $W/{\Bbb A}^1$,
 one has the genus function
 $$
  g(\eta)\;
  :=\;  r+1-|V(\Gamma)|
         + \sum_{v\in V(\Gamma_1)\cup V(\Gamma_2)}\,g(v)
 $$
 and the $H$-degree function
 $$
  d(\eta)\;
   :=\;
     \sum_{v\in V(\Gamma_1)}\,b_{\Gamma_1}(v)\cdot H|_{Y_1}
      + \sum_{v\in V(\Gamma_2)}\, b_{\Gamma_2}(v)\cdot H|_{Y_2}\,.
 $$
Denote by $|\eta|$ the triple of integers $(g(\eta), k_1+k_2; d(\eta))$.

On local standard charts of ${\mathfrak M}({\mathfrak W},(g,k;d))$,
 admissible triples $\eta$ with $|\eta|=(g,k;d)$ are used to encode
 at the topological level how a stable morphism $f$ from a {\it connected}
 curve $C$ to the degenerate fibers of $W[n]/{\Bbb A}^1[n]$ is to be
 realized as the gluing $f_1\sqcup f_2$ of relative stable morphisms $f_i$
 from a {\it possibly not connected} sub-curve $C_i$ of $C$ to fibers of
 some $Y[n_i]/{\Bbb A}^{n_i}$, $i=1,\,2$.
Even more importantly, each such $\eta$ corresponds to
 a {\it Cartier divisor} $({\mathbf L}_{\eta},{\mathbf s}_{\eta})$
 on ${\mathfrak M}({\mathfrak W},(g,k;d))$.
The zero-locus of ${\mathbf s}_{\eta}$ gives a closed ``divisoral"
 substack ${\mathfrak M}({\mathfrak W}_0,\eta)$ of
 ${\mathfrak M}({\mathfrak W},(g,k;d))$ that constitutes
 a (union of) component(s) of the fiber
 ${\mathfrak M}({\mathfrak W}_0,(g,k;d))$ of
 ${\mathfrak M}({\mathfrak W},(g,k;d))$ over ${\mathbf 0}\in {\Bbb A}^1$.
See
 [Li1: Definition 4.7, Proposition 4.12, Proposition 4.13] and
 [Li2: Definition 3.3].
Moreover, for each such $\eta$ one has a gluing morphism
 $$
  \Phi_{\eta}:
   {\mathfrak M}({\mathfrak Y}_1^{\it rel},\Gamma_1)
    \times_{E^r} {\mathfrak M}({\mathfrak Y}_2^{\rel},\Gamma_2)
  \rightarrow {\mathfrak M}({\mathfrak W},(g,k;d))\,,
 $$
 which is finite \'{e}tale of pure degree $|\Eq(\eta)|$ to
  its image
   ${\mathfrak M}
    ({\mathfrak Y}_1^{\rel}\,\sqcup\,{\mathfrak Y}_2^{\rel},\eta)$
  in and topologically isomorphic to
  ${\mathfrak M}({\mathfrak W}_0,\eta)$.
 [Li1: Sec.\ 4.2] and [Li2: Sec.\ 3.2].
(Here $\Eq(\eta)$ is the set of permutations of the $r$-many roots in
 $\Gamma_1$ that leaves $\eta$ unchanged.)

The obstruction theory associated to the deformation problems related
 to these moduli stacks are studied in [Li2: Sec.\ 1 and Sec.\ 5].
A cohomological description of the deformations of the separate
 constituents of a stable map and the natural clutching morphisms
 that relate the various separate deformation of the constituents
 are given there.
The perfectness of the obstruction theory and hence virtual fundamental
 classes on the various stacks involved here:
  $[{\mathfrak M}({\mathfrak W},(g,k;d))]^{\virt}$,
  $[{\mathfrak M}({\mathfrak Y}_i^{\rel},\Gamma_i)]^{\virt}$, ..., 
 are proved and constructed in [Li2: Sec.\ 2].
(See also [G-V: Sec.\ 2.8 and Sec.\ 2.9].)
These are then used to define the Gromov-Witten invariants of
 $Y_1\cup_D Y_2$  with values in ${\Bbb Q}$ and the relative
 Gromov-Witten invariants of the relative pairs $(Y_i,D)$
 with values in some $H_{\ast}(D^r,{\Bbb Q})$.

Given $(g,k;d)$, let
 $\Omega_{(g,k;d)}$ be the set of admissible triples $\eta$ for the
  gluing $Y_1\cup _D Y_2$ such that $|\eta|=(g,k;d)$,
 $\overline{\Omega}_{(g,k;d)}$ be the set of equivalence classes in
  $\Omega_{(g,k;d)}$ from re-ordering of roots, and
 ${\mathbf m}(\eta)$ be the product of the weight of roots of
  $\Gamma_1$ in $\eta\in\overline{\Omega}_{(g,k;d)}$.
Then J.\ Li's degeneration formula of Gromov-Witten invariants
 in numerical form reads: ([Li2: Sec.\ 3 and Sec.\ 4])

\bigskip

\noindent
{\bf Fact 1.3 [J.\ Li's degeneration formula].}
([Li2: Theorem 3.15 and Corollary 3.16].)
{\it
   $$
    \Psi_{(g,k;d)}^{W_t}(\alpha(t),\zeta)\;
     =\; \sum_{\eta\in\overline{\Omega}_{(g,k;d)}}\;
        \frac{{\mathbf m}(\eta)}{|\Eq(\eta)|}\,
          \sum_{j\in K_{\eta}}\,
           \left[ \Psi_{\Gamma_1}^{Y_1^{\rel}}(j_1^{\ast}\alpha(0),
                    \zeta_{\eta,1,j})\,
                  \bullet\,
                  \Psi_{\Gamma_2}^{Y_2^{\rel}}(j_2^{\ast}\alpha(0),
                    \zeta_{\eta,2,j})
           \right]_0\,.
  $$
} % end-fact

\bigskip

\noindent
Here the various $\Psi$'s are the Gromov-Witten invariants
 (resp.\ relative Gromov-Witten classes) in ${\Bbb Q}$
 (resp.\ in various $H_{\ast}(D^r,{\Bbb Q})$'s)
 that arise from intersections with the virtual fundamental classes
 constructed in [Li2: Sec.\ 2] of cycles from pull-back via
 evaluation maps or projection maps of cycles on $W/{\Bbb A}^1$,
 $Y_1$, $Y_2$ or on the related moduli stack of stable curves.
See [Li2: Sec.\ 0] for notations unexplained here with
   $(g,k;d)$ here = $\Gamma$ there,
   $\Psi^{W_t}_{(g,k;d)}$ here = $\Phi^{W_t}_{\Gamma}$ there, and
   $\zeta$ here = $\beta$ there
 to avoid confusions with our later use of notations.
(See also Corollary 2.2 in Sec.\ 2
 % Corollary ??? [J Li's degeneration formula']
 for a nearly complete statement of the formula in our case.)

We only want to remark that $\Psi^{W_t}_{(g,k;d)}$ is the {\it sum} of
 the usual Gromov-Witten invariants defined via
 $\overline{\cal M}_{g,k}(W_t,\beta)$
 over $\beta\in A_1(W_t)$ (or $\in H_2(W_t;{\Bbb Z})$) such that
 $\beta\cdot H|_{W_t}=d$.
As $H$ varies, the degeneration formula above sums over a {\it different}
 collection of curve classes in the Mori cone $\overline{\NE}(W_t)$.
For other new applications of J.\ Li's formula, one needs
 a form of the degeneration formula not with respect to
 $H$-degree $d$ but with respect to a curve class $\beta$ in $A_1(W_t)$
 or $H_2(W_t;{\Bbb Z})$.
Due to the possible monodromy effect, we do not know as yet
 what such a formula should look like in the most general form. 
But for the degeneration $W\rightarrow {\Bbb A}^1$ that comes
 from blowing up a trivial family $X\times{\Bbb A}^1$, J.\ Li's formula
 can be readily modified into the desired form once one traces through
 the details of [Li2; Sec.\ 3 and Sec.\ 4].
We now explain this modification.

\bigskip

\section{A degeneration formula with respect to a curve class.}
Given a smooth projective variety $X$ and a smooth subvariety
 $Z\subset X$, let $p:W\rightarrow X\times {\Bbb A}^1$
 be the blow-up of $X\times{\Bbb A}^1$ along $Z\times{\mathbf 0}$,
 where ${\mathbf 0}$ is a point on ${\Bbb A}^1$.
The induced map $\pi: W\rightarrow {\Bbb A}^1$ gives a degeneration
 of the kind discussed in [Li1] with $W_0=Y_1\cup_E Y_2$,
 where
  $Y_1=\Bl_Z X$,
  $Y_2={\Bbb P}({\cal N}_{Z/X}\oplus{\cal O}_Z)$, and
  $E={\Bbb P}{\cal N}_{Z/X}$
  with ${\cal N}_{Z/X}$ being the normal bundle of $Z$ in $X$.
We denote also the induced morphisms as
 $p: W_0\rightarrow X$, $p_1: Y_1\rightarrow X$, and
 $p_2: Y_2\rightarrow X$.
Since $X$ is projective, there exists a relative ample line bundle $H$
 on $W/{\Bbb A}^1$, which will be fixed in this section.
$X$ will be identified as a general fiber of $W/{\Bbb A}^1$ as well
 whenever necessary.
For simplicity of notations, we assume that $Z$ is connected.

Given a triple $(g,k;d)$, let
 $C_{(H,\,d)}(X)
   =\{\,\beta \in A_1(X)\,:\, H\cdot \beta =d\,\}$.
Then
 $$
  {\mathfrak M}({\mathfrak W},(g,k;d))\;
   =\; \coprod_{\beta\in C_{(H,\,d)}(X)}
              {\mathfrak M}({\mathfrak W}, (g,k;\beta) )\,,
 $$
 where ${\mathfrak M}({\mathfrak W},(g,k;\beta))$ is the stack of
  stable morphisms from prestable curves of genus $g$ with $k$ marked
  points to the universal family of the stack ${\mathfrak W}$ of
  expanded degenerations associated to $W/{\Bbb A}^1$ such that
  after the post-composition with the morphisms
  ${\mathfrak W}\rightarrow X\times{\Bbb A}^1\rightarrow X$,
  the images of the stable morphisms lie in the curve class $\beta$.
Each ${\mathfrak M}({\mathfrak W},(g,k;\beta))$ is a Deligne-Mumford
 stack as in [Li1: Theorem 3.10].
Being a union of connected components of
 ${\mathfrak M}({\mathfrak W}, (g,k;d))$,
 the deformation-obstruction theory on
 ${\mathfrak M}({\mathfrak W},(g,k;\beta))$ and
 ${\mathfrak M}({\mathfrak W}_0,(g,k;\beta))
   := {\mathfrak M}({\mathfrak W},(g,k;\beta))
               \times_{{\scriptsizeBbb A}^1}{\mathbf 0}$
 are the same as those in [Li2: Sec.\ 1].
The tangent-obstruction complex on
  ${\mathfrak M}({\mathfrak W},(g,k;\beta))$,
 its perfectness and its resolving locally free $2$-term complex
  ${\mathbf E}^{\bullet}=[E^1\rightarrow E^2]$,
 the relative Kuranishi structure on ${\mathbf E}^{\bullet}$,
 the associated cone class in $E^2$, and
 the virtual fundamental class
  $[{\mathfrak M}({\mathfrak W},(g,k;\beta))]^{\virt}$
  defined via the cone class and Gysin pull-back via the $0$-section
  of $E^2$ over ${\mathfrak M}({\mathfrak W},(g,k;\beta))$
 are those constructed on [Li2: Sec.\ 1 and Sec.\ 2] for
 ${\mathfrak M}({\mathfrak W},(g,k;d))$ but restricted to its connected
 components that constitute ${\mathfrak M}({\mathfrak W},(g,k;\beta))$.

Recall
 the set $\Omega_{(g,k;d)}$ of admissible triples $\eta$ such that
  $|\eta|=(g,k;d)$ and
 the quotient set $\overline{\Omega}_{(g,k;d)}$ reviewed in Sec.\ 1.
For an admissible weighted graph $\Gamma$ for a relative pair,
 define $b(\Gamma):=\sum_{v\in V(\Gamma)} b(v)$.
Define the {\it $\beta$-compatible subset} of $\Omega_{(g,k;d)}$ by
 $$
  \Omega_{(g,k;\beta)}^H\;
   :=\; \left\{\,
    \eta=(\Gamma_1,\Gamma_2,I)\in\Omega_{(g,k;d)}\,\left|\,
       p_{1\ast}b(\Gamma_1)+p_{2\ast}b(\Gamma_2)=\beta
         \right.
      \,\right\}\,.
 $$
To proceed, we need to understand the Cartier divisors
 $({\mathbf L},{\mathbf s}_{\eta})$ associated to
 $\eta\in\Omega_{(g,k;d)}$.

Let $\eta=(\Gamma_1,\Gamma_2,I)$.
Recall the substack
 ${\mathfrak M}({\mathfrak Y}_1^{\rel}\,\sqcup\,{\mathfrak Y}_2,\eta)$
 in ${\mathfrak M}({\mathfrak W},(g,k;d))$, which is the image of
 the gluing morphism $\Phi_{\eta}$ associated to $\eta$.
For any local chart $S$ of ${\mathfrak M}({\mathfrak W},(g,k;d))$
(i.e.\ a morphism $S\rightarrow {\mathfrak M}({\mathfrak W},(g,k;d))$,
 whose associated universal family will be denoted by $f$),
 define
 $$
  S_{\eta}\;
   := \; S \times_{{\mathfrak M}({\mathfrak W},(g,k;d))}\,
         {\mathfrak M}
         ({\mathfrak Y}_1^{\rel}\,\sqcup\,{\mathfrak Y}_2^{\rel},\eta)\,.
 $$
(For readers not familiar with stacks, it is instructive to think of
 $S_{\eta}$ as the intersection of $S$ with
 ${\mathfrak M}
  ({\mathfrak Y}_1^{\rel}\,\sqcup\,{\mathfrak Y}_2^{\rel},\eta)$
 in ${\mathfrak M}({\mathfrak W},(g,k;d))$,
 $S_{\eta}$ is an \'{e}tale subscheme of $S$.)
For $S$ with $S_{\eta}$ empty, $({\mathbf L}_{\eta},{\mathbf s}_{\eta})|_S$
 is defined to be the trivial Cartier divisor $({\cal O}_S,1)$.
In general, $S$ can be covered by local \'{e}tale charts $S_{\alpha}$
 such that the restriction/pull-back of $f$ on each $S_{\alpha}$ is
 represented by a square
 $f_{\alpha}: {\cal C}_{\alpha}/S_{\alpha}
              \rightarrow W[n_{\alpha}]/{\Bbb A}^1[n_{\alpha}]$
 in such a way that
 \begin{itemize}
  \item[(1)]
   The induced map $S_{\alpha\eta}\rightarrow {\Bbb A}^1[n_{\alpha}]$
    factors through $S_{\alpha\eta}\rightarrow {\mathbf H}_{l_{\alpha}}$,
    where ${\mathbf H}_{l_{\alpha}}$ is a coordinate hyperplane in
    ${\Bbb A}^1[n]={\Bbb A}^{n+1}$.

  \item[(2)]
   Recall the distinguished locus
    ${\mathbf D}_{l_{\alpha}}\simeq {\mathbf H}_{l_{\alpha}}\times D$
    over ${\mathbf H}_{l_{\alpha}}$ in $W[n_{\alpha}]$
    ([Li1: Sec.\ 1.1]).
   Then $f_{\alpha}^{-1}({\mathbf D}_{l_{\alpha}})$ divides the
    (connected) prestable curve ${\cal C}_{\alpha}/S_{\alpha\eta}$
    into two collections of (possibly not connected) prestable curves
     ${\cal C}_{\alpha,1}/S_{\alpha\eta}$ and
     ${\cal C}_{\alpha,2}/S_{\alpha\eta}$,
    realizing ${\cal C}_{\alpha}/S_{\alpha\eta}$ as gluing of
    a curve of type $\Gamma_1$ and a curve of type $\Gamma_2$
    along the distinguished divisors 
     $f_{\alpha}^{-1}({\mathbf D}_{l_{\alpha}})$ now on
     ${\cal C}_{\alpha,1}/S_{\alpha\eta}$ and
     ${\cal C}_{\alpha,2}/S_{\alpha\eta}$ respectively.
 \end{itemize}
(See [Li1: Sec.\ 1.1 and Sec.\ 2.2] or {\sc Figure} 2-1
  for the locus ${\mathbf D}_l\subset W[n]$
  over ${\mathbf H}_l \subset {\Bbb A}^1[n]$
  - here we retain J.\ Li's notation ${\mathbf D}_l$ for easy referring
    but, to be consistent with our notations, we should denote it
    by ${\mathbf E}_l$.)
Let $(L_{l_\alpha},s_{l_{\alpha}})$ be (any) Cartier divisor on
 ${\Bbb A}^1[n_{\alpha}]$ with the zero-scheme of $s_{l_{\alpha}}$
 being the hyperplane ${\mathbf H}_{l_{\alpha}}$.
Then $f_{\alpha}^{\ast}(L_{l_{\alpha}},s_{l_{\alpha}})$
 is a Cartier divisor
 $({\mathbf L}_{\eta,\alpha},{\mathbf s}_{\eta,\alpha})$ on $S_{\alpha}$
 with zero-locus $S_{\alpha\eta}$.
Let $S$ run over the charts in an atlas of
 ${\mathfrak M}({\mathfrak W},(g,k;d))$,
 one obtains a collection of local Cartier divisors.
They  descend to a Cartier divisor
 $({\mathbf L}_{\eta},{\mathbf s}_{\eta})$
 on the stack ${\mathfrak M}({\mathfrak W},(g,k;d))$,
 [Li2: Sec.\ 3.1 and Lemma 3.4].

%\marginpar{\raggedright\tiny $\bullet$
% {\sc Figure} ???: \newline
%  {\tt exp-dgn.eps}}

\begin{figure}[htbp]
 \setcaption{{\sc Figure} 2-1.
  \baselineskip 14pt
  A local model $W[n]/{\Bbb A}^1[n]$ for the universal family of
   the stack ${\mathfrak W}$ of expanded degenerations associated
   to a degeneration $W/{\Bbb A}^1$ and
  its tautological morphism to $W/{\Bbb A}^1$.
  Over each coordinate hyperplane ${\mathbf H}_l$ in the base
   ${\Bbb A}^1[n]:= {\Bbb A}^{n+1}$ is a distinguished locus
   ${\mathbf D}_l$ in $W[n]$, crossing the singular locus of the fibers
   of $W[n]/{\Bbb A}^1[n]$ and constituting a connected component of
   the subscheme of such loci in $W[n]$.
  Illustrated here is the case $n=1$.
 } % end-setcaption
\centerline{\psfig{figure=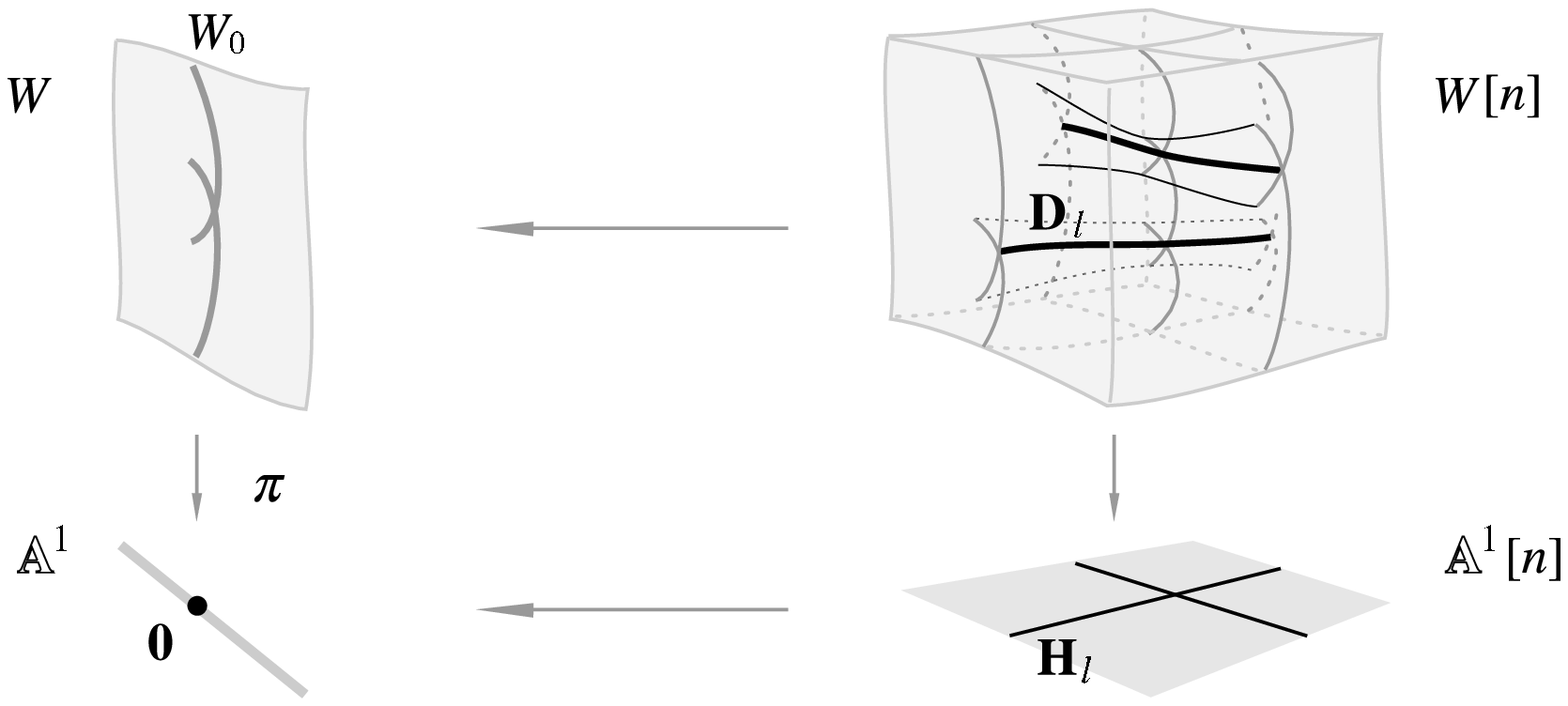,width=13cm,caption=}}
\end{figure}

Recall also the Cartier divisor $({\mathbf L}_0,{\mathbf t}_0)$
 on ${\mathfrak M}({\mathfrak W},(g,k;d))$ that is the pull-back
 of the Cartier divisor $({\cal O}_{{\scriptsizeBbb A}^1},t)$
 on ${\Bbb A}^1$,
   where ${\Bbb A}^1=\Spec {\Bbb C}[t]$ and
         ${\mathbf 0}=(t)\in \Spec{\Bbb C}[t]$,
 via the tautological morphism
 ${\mathfrak M}({\mathfrak W},(g,k;d))\rightarrow {\Bbb A}^1$.
Then: 

\bigskip

\noindent
{\bf Lemma 2.1 [validity].} {\it
 J.\ Li's results on ${\mathfrak M}({\mathfrak W},(g,k;d))\,:$
  Proposition 3.5, Theorem 3.6, Lemma 3.10 - 3.14, Theorem 3.15, and
  Corollary 3.16 in {\rm [Li2: Sec.\ 3]}
 hold for the substack ${\mathfrak M}({\mathfrak W},(g,k;\beta))$
 as well, with $\eta$ in the quote running over
 $\overline{\Omega}_{(g,k;\beta)}^H$.
} % end-lemma

\bigskip

\noindent
{\it Proof.}
 From the above highlight of J.\ Li's construction of
  $({\mathbf L}_{\eta},{\mathbf s}_{\eta})$,
 it is immediate that 
  for $\eta\in \Omega_{(g,k;d)}-\Omega_{(g,k;\beta)}^H$,
 the restriction of $({\mathbf L}_{\eta},{\mathbf s}_{\eta})$ to
  ${\mathfrak M}({\mathfrak W},(g,k;\beta))$ is the trivial
  Cartier divisor
  $({\cal O}_{{\mathfrak M}({\mathfrak W},(g,k;\beta))},{\mathbf 1})$.
 Similarly, for $\eta\in\Omega_{(g,k;\beta)}^H$, the restriction
  of $({\mathbf L}_{\eta},{\mathbf s}_{\eta})$ to
  $\mathfrak M({\mathfrak W},(g,k;d))
               -{\mathfrak M}({\mathfrak W},(g,k;\beta))$
  is a trivial Cartier divisor as well.
 This implies that J.\ Li's results can be applied to
  ${\mathfrak M}({\mathfrak W},(g,k;\beta))$ and
   ${\mathfrak M}({\mathfrak W}, (g,k;d))
                -{\mathfrak M}({\mathfrak W},(g,k;\beta))$
  respectively alone.
 The lemma then follows.

\noindent\hspace{14cm}$\Box$

\bigskip

Explicitly,
$\otimes_{\eta\in\overline{\Omega}_{(g,k;\beta)}^H}
  ({\mathbf L}_{\eta},{\mathbf s}_{\eta})\;
  = ({\mathbf L}_0,{\mathbf t}_0)$
on ${\mathfrak M}({\mathfrak W},(g,k;\beta))$, cf.\ [Li2: Proposition 3.5].
This implies the following version of the degeneration formula,
 cf.\ [Li2: Theorem 3.6] and notations in Corollary 2.2 below:

{\small
$$
 \Psi_{(g,k;\beta)}^X(\alpha,\zeta)\;
  =\; {\mathbf q}_{\ast 0}\left(
        \sum_{\eta\in\overline{\Omega}_{(g,k;\beta)}^H}\,
           \left( \rule{0ex}{3ex}
             \ev_0^{\ast}\left(\alpha(0)\right)
              \cdot \pi_{g,k}^{\ast}(\zeta)
              \cdot \left(\rule{0ex}{2.4ex}
                     c_1({\mathbf L}_{\eta},{\mathbf s}_{\eta})
                      [{\mathfrak M}({\mathfrak W},(g,k;\beta))]^{\virt}
                    \right)
           \right)
                          \right)\,;
$$
} % end-small

\noindent
see [Li2: Sec.\ 3.1] for more explanations of notations.
The obstruction theory on ${\mathfrak M}({\mathfrak W}_0,\eta)$
 is constructed in [Li2: Sec.\ 3.2] and
the following identities hold for
 ${\mathfrak M}({\mathfrak W}_0,(g,k;\beta))$,
cf.\ [Li2: Sec.\ 3.2, Lemma 3.10-14, Theorem 3.15]:
\begin{eqnarray*}
 c_1({\mathbf L}_0,{\mathbf t}_0)
   [{\mathfrak M}({\mathfrak W},(g,k;\beta))]^{\virt}
  & = & [{\mathfrak M}({\mathfrak W}_0,(g,k;\beta))]^{\virt}\,, \\[.6ex]
 c_1({\mathbf L}_{\eta},{\mathbf s}_{\eta})
      [{\mathfrak M}({\mathfrak W},(g,k;\beta))]^{\virt}
 & = & [{\mathfrak M}({\mathfrak W}_0,\eta)]^{\virt}\;
   =\; {\mathbf m}(\eta)\,
       [{\mathfrak M}({\mathfrak Y}_1^{\rel}\,
             \sqcup\,{\mathfrak Y}_2^{\rel},\eta)]^{\virt}\,,   \\[.6ex]
 [{\mathfrak M}({\mathfrak W}_0,(g,k;\beta))]^{\virt}
  & = & \sum_{\eta\in\overline{\Omega}_{(g,k;\beta)}^H}\,
         {\mathbf m}(\eta)\,
         [{\mathfrak M}({\mathfrak Y}_1^{\rel}\,
                         \sqcup\,{\mathfrak Y}_2^{\rel},\eta)]^{\virt}\,,
\end{eqnarray*}
and
$$
 \frac{1}{|\Eq(\eta)|}\Phi_{\eta\ast}\Delta^!
  \left(
   [{\mathfrak M}({\mathfrak Y}_1^{\rel, \Gamma_1})]^{\virt}
    \times [{\mathfrak M}({\mathfrak Y}_2^{\rel},\Gamma_2)]^{\vert}
  \right)\;
  =\; [{\mathfrak M}
       ({\mathfrak Y}_1^{\rel}\,\sqcup\,{\mathfrak Y}_2^{\rel},
                                                  \eta)]^{\virt}\,.
$$
These together imply a degeneration formula for $\Psi_{(g,k;\beta)}^X$
 in cycle form, which can be converted into the equivalent numerical form,
cf.\ [Li2: Sec.\ 3.2 and Sec.\ 4].

To summarize, given an admissible triple
 $\eta=(\Gamma_1,\Gamma_2,I)\in\Omega_{(g,k;\beta)}^H$ with $r$-many
  roots in $\Gamma_i$,
 recall the evaluation morphism from distinguished marked points:
  ${\mathbf q}_i:{\mathfrak M}({\mathfrak W},\Gamma_i)\rightarrow E^r$,
 [Li2: Sec.\ 0].
Then [Li2: Sec.\ 3 and Sec.\ 4] and Lemma 2.1 above imply that:
 (cf.\ [Li2: Theorem 3.15 and Corollary 3.16])

\bigskip

\noindent
{\bf Corollary 2.2 [J.\ Li's degeneration formula$^\prime$].}
{\it
 Fix a relative ample line bundle $H$ on $W$.
 Let
  $\beta\in A_1(X)$,
  $d :=\beta\cdot H$,
  $\Omega_{(g,k;\beta)}^H$ be the subset of admissible triples
   $(\Gamma_1,\Gamma_2,I)$ in $\Omega_{(g,k;d)}$ such that
   $p_{1\ast}b(\Gamma_1)+p_{2\ast}b(\Gamma_2)=\beta$,
  $\alpha\in A_{\ast}(X)^{\times k}$,
   whose extension to $H^0({\Bbb A}^1;R\pi_{\ast}{\Bbb Q}_W)^{\times k}$
   is denoted still by $\alpha$, and
  $\zeta\in A_{\ast}({\mathfrak M}_{g,k})$.
 Denote by $\Psi_{(g,k;\beta)}^X(\alpha,\zeta)$ 
  the usual Gromov-Witten invariant of $X$ associated to these data.
 For $\eta\in\overline{\Omega}_{(g,k;\beta)}^H$, assume that
  $G_{\eta}^{\ast}(\zeta)
    =\sum_{j\in K_{\eta}}\zeta_{\eta,1,j}\bboxtimes \zeta_{\eta,2,j}$,
  where
   $G_{\eta}: {\mathfrak M}_{\Gamma_1^o}\times {\mathfrak M}_{\Gamma_2^o}
              \rightarrow {\mathfrak M}_{g,k}$ is the natural morphism
   between the related moduli stack of nodal curves
   $($cf.\ {\rm [Li2: Sec.\ 0]} for more explanations$)$.
 Then
 $$
  \Psi_{(g,k;\beta)}^X(\alpha,\zeta)\;
   =\; \sum_{\eta\in\overline{\Omega}_{(g,k;\beta)}^H}\;
        \frac{{\mathbf m}(\eta)}{|\Eq(\eta)|}\,
          \sum_{j\in K_{\eta}}\,
           \left[ \Psi_{\Gamma_1}^{Y_1^{\rel}}(j_1^{\ast}\alpha(0),
                    \zeta_{\eta,1,j})\,
                  \bullet\,
                  \Psi_{\Gamma_2}^{Y_2^{\rel}}(j_2^{\ast}\alpha(0),
                    \zeta_{\eta,2,j})
           \right]_0\,,
 $$
 where
  $j_i:Y_i\rightarrow W_0$ is the inclusion map,
  $$
   \Psi_{\Gamma_i}^{Y_i^{rel}}(j_i^{\alpha}(0),\zeta_{\eta,i,j})
   = {\mathbf q}_{i\,\ast}\left(
       \ev^{\ast}(j_i^{\ast}\alpha(0))
      \cdot \pi_{\Gamma_i}^{\ast}(\zeta_{\eta,i,j})
      \cdot  [{\mathfrak M}({\mathfrak Y}_i^{\rel},\Gamma_i)]^{\virt}
                          \right) \in H_{\ast}(E^r)\;, \;\; i=1,2\,,
  $$
  {\rm (}here
   $\pi_{\Gamma_i}:{\mathfrak M}({\mathfrak Y}_i^{\rel},\Gamma_i)
       \rightarrow {\mathfrak M}_{\Gamma_i^o}$,
   {\rm [Li2: Sec.\ 0]}{\rm )},
  $\bullet$ is the intersection product on the relevant 
    $A_{\ast}(E^r)_{\scriptsizeBbb Q}$ for each summand, and
  $[\,\cdot\,]_0$ is the degree-$0$ component of a cycle class.
 
 In cycle form,
 $$
  [{\mathfrak M}({\mathfrak W}_0,(g,k;\beta))]^{\virt}\;
   =\; \sum_{\eta\in\overline{\Omega}_{(g,k;\beta)}^H}\;
         \frac{{\mathbf m}(\eta)}{|\Eq(\eta)|}\,
             \Phi_{\eta\ast}\Delta^!
               \left(
                [{\mathfrak M}({\mathfrak Y}_1^{\rel},\Gamma_1)]^{\virt}
                \times
                [{\mathfrak M}({\mathfrak Y}_2^{\rel},\Gamma_2)]^{\virt}
               \right)\,,
 $$
 where
 $\Delta^!$ is the Gysin map associated to the diagonal map
 $\Delta:E^r\rightarrow E^r\times E^r$ for the relevant $E^r$
  in each summand.
} % end-lemma

\bigskip

\noindent
With $A_1(X)$ replaced by $H_2(X;{\Bbb Z})$ in the discussion,
 identical degeneration formulas with respect to
 $\beta\in H_2(X;{\Bbb Z})$ follow.

In these identities, though $\beta$ and hence the left-hand side are
 now independent of $H$, $\Omega_{(g,k;\beta)}^H$ and hence
 the right-hand side may still depend on $H$ since it is chosen as
 a subset of $\Omega_{(g,k;d)}$ and the latter does depend on $H$.
Let us now turn to this issue in Sec.\ 3.

\bigskip

\section{The $H$-(in)dependence of $\Omega_{(g,k;\beta)}^H$.}

In this section we discuss how to stabilize $\Omega_{(g,k;\beta)}^H$
 that appears in the degeneration formula in Corollary 2.2
 % Corollary ??? [J.\ Li's degeneration formula']
 and potentially depends on the choice of a relative ample line bundle
 $H$ on $W/{\Bbb A}^1$.

With the same notation as in Sec.\ 2,
 every line bundle $L$ on $X$ induces canonically
 a line bundle
 $\widehat{L}$ on $X\times{\Bbb A}^1$ by pulling back.
Since ${\Bbb A}^1$ is affine, whose coordinate ring ${\Bbb C}[t]$
 separates points and tangent vectors, very-ampleness of $L$ on $X$
 implies very-ampleness of $\widehat{L}$ on $X\times{\Bbb A}^1$.
Since $W$ is a blow-up of $X\times{\Bbb A}^1$ with exceptional
 divisor $Y_2$, it follows from [Ha] or [C-H] that
 if $L$ is sufficiently very ample on $X$, then
 $H := p^{\ast}\widehat{L}\otimes {\cal O}_W(-Y_2)$
 is very ample on $W$ (and hence $\pi$-ample).

\bigskip

\noindent
{\it Remark 3.1 {\rm [}sufficiently very ampleness{\rm ]}.}
 Let $L^{\prime}$ be a very ample line bundle on $X$.
 Associated to $L^{\prime}$ is a homogeneous coordinate ring $R$
  for $X\times{\Bbb A}^1$ associated to $\widehat{L^{\prime}}$.
 (Or one may consider the compactification $X\times{\Bbb P}^1$
  of $X\times{\Bbb A}^1$
  and define $\widehat{L^{\prime}}$ to be
  $L^{\prime}\bboxtimes {\cal O}_{{\scriptsizeBbb P}^1}(1)$.).
 Suppose that the ideal sheaf ${\cal I}$ of $Z\times{\mathbf 0}$
  in $W$ is (finitely) generated by homogeneous elements of $R$
  of degree $\le a$, then [C-H] says that
  $p^{\ast}\widehat{L^{\prime}}^{\otimes c}\otimes{\cal O}_W(-Y_2)$
  is very ample for $c\ge a+1$.
 Take now $L=L^{\prime\,\otimes c}$ with $c\ge a+1$.

\bigskip

\noindent
{\bf Lemma 3.2 [$H$-independence of $\Omega_{(g,k;\beta)}^H$].}
{\it
 Let $H$ be a relative ample line bundle on $W/{\Bbb A}^1$ that is
  associated to a sufficiently very ample line bundle $L$ on $X$
  as above.
 Then the set $\overline{\Omega}_{(g,k;\beta)}^H$ that appears
  in the degeneration formula in Corollary 2.2
  % Corollary ??? [J.\ Li's degeneration formula']
  does not depend on $L$ or $H$.
 Rather, it depends only on $\beta\in A_1(X)$
  and the normal bundle ${\cal N}_{Z/X}$ of $Z$ in $X$.
} %

\bigskip

\noindent
{\it Proof.}
Fix such an $L$ on $X$, then the associated $H$ is very ample on $W$.
Let
 $H_1 := H|_{Y_1}$ and
 $H_2 := H|_{Y_2}$.
Recall the morphisms
 $p_1: (Y_1,E)\rightarrow (X,Z)$,
 $p_2: (Y_2,E)\rightarrow (X,Z)$, and
 $p: Y_1\cup_E Y_2\rightarrow X$.
We assume that $Z$ has codimension $r$ in $X$.
To simplify notations, we shall not distinguish a line bundle and
 its associated equivalence class of Cartier divisors in our discussion.
Denote by $\gamma$ the line class in a ${\Bbb P}^{\,r-1}$-fiber of $E/Z$.
Then
 both $\NE(p_1)$ and $\NE(p_2)$ are generated by $\gamma$ and
 ${\Bbb R}_+\gamma$ is an extremal ray in both Mori cones
  $\overline{\NE}(Y_1)$ and $\overline{\NE}(Y_2)$.
(Caution that $E\cdot \gamma=-1$ as cycles in $Y_1$
  while $E\cdot\gamma =+1$ as cycles in $Y_2$.)

%\marginpar{\raggedright\tiny $\bullet$
% {\sc Figure} ???: \newline
%  {\tt moricone.eps}}

\begin{figure}[htbp]
 \setcaption{{\sc Figure} 3-1.
  \baselineskip 14pt
   The Mori cones $\overline{\NE}(X)$, $\overline{\NE}(Y_1)$,
    and $\overline{\NE}(Y_2)$ in the problem,
   the morphisms between them, and
   the extremal ray ${\Bbb R}_+\gamma$ associated to
    $p_i:Y_1\rightarrow X$, $i=1,\,2$,
   are illustrated.
  Both the cross-section of the hyperplane
   $\{H_i=d_i\}$ in $N_1(Y_i)_{\scriptsizeBbb R}$ with
   $\overline{\NE}(Y_i)$,
   where $d_i :=H_i\cdot \widetilde{\beta}_i$, $i=1,\,2$, and
  the ray $p_{i\ast}^{-1}(p_{i\ast}\widetilde{\beta}_i)$
   are indicated.
 } % end-setcaption
\centerline{\psfig{figure=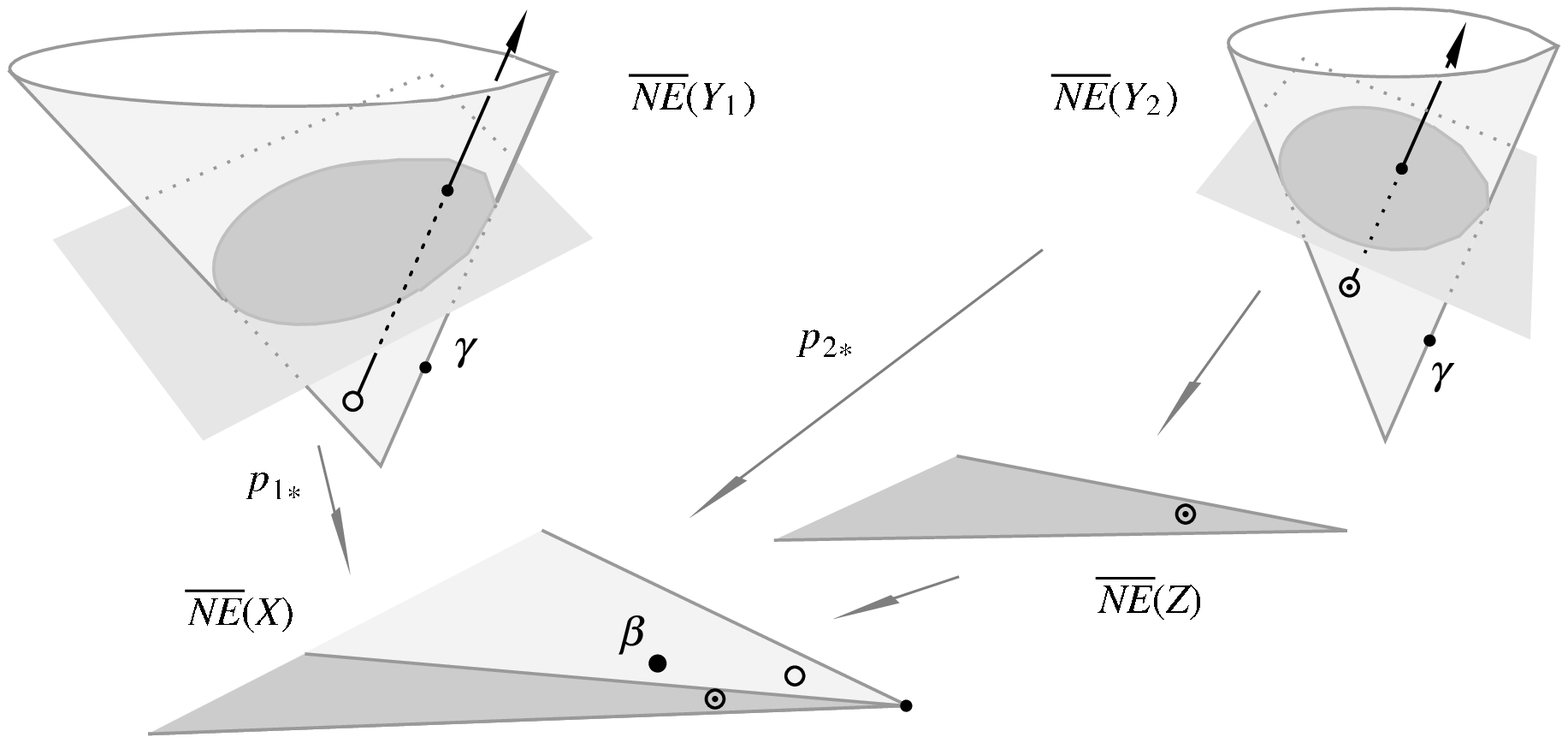,width=13cm,caption=}}
\end{figure}

By construction,
 $$
  H_1 \;
   =\; p_1^{\ast}L \otimes{\cal O}_{Y_1}(-E)\;
    =\; p_1^{\ast}L -E\,.
 $$
For $H_2$, observe that:
 $H_2\cdot \gamma$ on $Y_2$ $= H\cdot\gamma$ on $W$
 $= -Y_2\cdot\gamma$ on $W$ $= +1$,
 where
  we have used the projection formula for push-pull of cycles under
   a proper morphism and
  the fact that $Y_2\cdot \gamma$ on $W$ is the degree of
   ${\cal O}_W(Y_2)|_{\gamma}$ on $\gamma$
   and the latter bundle is isomorphic to 
   ${\cal O}_{Y_2/Z}(-1)|_{\gamma}
                           \simeq {\cal O}_{{\scriptsizeBbb P}^1}(-1)$.
Together with the structure theorem of the Chow ring of projective space
 bundles, it follows that
 $$
  H_2 \;
   =\; p_2^{\ast}(L\cdot Z) + c_1({\cal O}_{Y_2/Z}(1))\;
   =\; p_2^{\ast}(L\cdot Z) + E\,.
 $$
$H_1$ and $H_2$ are very ample on $Y_1$ and $Y_2$ respectively.

Let $\eta=(\Gamma_1,\Gamma_2,I)\in \Omega_{(g,k;\beta)}^H$.
Then the pairs
 $(b(\Gamma_1),b(\Gamma_2))$ are characterized by the conditions:
 $$
  p_{1\ast} b(\Gamma_1)+p_{2\ast} b(\Gamma_2)\;
  =\;\beta
  \hspace{1em}\mbox{and}\hspace{1em}
  H_1\cdot b(\Gamma_1)+H_2\cdot b(\Gamma_2)\;
  =\; L\cdot \beta\; =\; d\,.
 $$
Since when $\beta=0$ or $\notin\NE(X)$, the statement in the Lemma
 holds vacuously, we shall assume that $0\ne \beta\in\NE(X)$ and
 that $b(\Gamma_i)\in \NE(Y_i)$, $i=1,\,2$, in the following discussion.

Let
 $(\widetilde{\beta}_1,\widetilde{\beta}_2)
    \in \NE(Y_1)_{\scriptsizeBbb Z}\times\NE(Y_2)_{\scriptsizeBbb Z}$
 be a pair of curves that satisfies the curve condition
 $p_{1\ast}\widetilde{\beta}_1+p_{2\ast}\widetilde{\beta}_2=\beta$.
Consider the following cases.

\bigskip

\noindent 
{\it Case $(a):$
 $\widetilde{\beta}_1\ne 0\;$ and $\;\widetilde{\beta}_2\ne 0$.}
Then
 $p_{i\ast}^{-1}(p_{i\ast}\widetilde{\beta_i}) \cap\NE(Y_i)$
 is a ray parallel to the extremal ray ${\Bbb R}_+\gamma$ of
 $\overline{\NE}(Y_i)$.
Let $\widetilde{\beta}_i^{\,0}$ be the first (non-apex) lattice point
 on this ray, then
 $\widetilde{\beta}_i = \widetilde{\beta}_i^{\,0} + l_i \gamma$,
 for a unique $l_i\ge 0$, $i=1,\,2$, cf.\ {\sc Figure} 3-1.
 % Figure ??? : moricone.eps
The $H$-degree condition implies that
 $(\widetilde{\beta}_1,\widetilde{\beta}_2)=(b(\Gamma_1),b(\Gamma_2))$
 for a $(\Gamma_1,\Gamma_2,I)\in\Omega_{(g,k;\beta)}^H$ if and only if
 $$
  E\cdot \widetilde{\beta}_1^{\,0}\,-\, l_1\;
  =\;  l_2 \;
  =\;  \mu_{\eta}\,
  :=\, \mbox{total root weight $\sum_i\mu_{1,i}$ of $\Gamma_1$}\,.
 $$
These conditions on the pairs $(l_1,l_2)$ of non-negative integers
 are independent of $L$ and hence $H$.

Since the weight functions
 $b:\Gamma_i\rightarrow A_1(Y_i)/\sim_{\scriptsizealg}$, $i=1,\,2$,
 are subject only to the conditions $b(\Gamma_i)=\widetilde{\beta}_i$,
 the set of their possible choices is determined only by
 $(\widetilde{\beta}_1,\widetilde{\beta}_2)$ and the semi-groups:
 $\NE(Y_i)_{\scriptsizeBbb Z}$, $i=1,\,2$.
Again this is independent of $L$ and hence $H$.
The choices of weighted roots in $\Gamma_1$ and $\Gamma_2$ are only
 constrained by the condition:
 $l_2=$ their common total weight $\mu_{\eta}$.
This is irrelevant to $H$ as well.
Finally, the choices of legs and $I$ are irrelevant to $H$ by nature.
This concludes the discussion for Case (a). 

\bigskip

\noindent
{\it Case $(b):$ $\widetilde{\beta}_2=0$ or $\widetilde{\beta}_1=0$.}
In the first situation, apply the above discussion for $Y_1$,
 one concludes that
 $\widetilde{\beta}_1
  =\widetilde{\beta}_1^{\,0}+ (E\cdot\widetilde{\beta}_1^{\,0})\gamma$.
In the second situation, $\beta$ is contained in $Z$ and
 $\widetilde{\beta}_2 = \widetilde{\beta}_2^{\,0}$.
Either way, the same conclusion as in Case (a) follows.
This concludes the proof.

\noindent\hspace{14cm}$\Box$

\bigskip

Denote the stabilized $\Omega_{(g,k;\beta)}^H$ by $\Omega_{(g,k;\beta)}$
 and its set of equivalence classes by $\overline{\Omega}_{(g,k;\beta)}$.
Denote with a superscript $^{0}$ for the first (non-cone-apex) lattice
 point of a ray in $NE(Y_i)$ parallel to ${\Bbb R}_+\gamma$.
Then the proof of Lemma 3.2 characterizes $\Omega_{(g,k;\beta)}$
 for $0\ne \beta\in NE(X)$ as:
(The set $\Omega_{(g,k;\,0)}$ is immediate.)

{\footnotesize
 $$
  \Omega_{(g,k;\beta)}\;
  =\; \coprod_{\tiny
        \begin{array}{c}
          (\widetilde{\beta}_1^{\,0},\,\widetilde{\beta}_2^{\,0})
                \in NE(Y_1)\times NE(Y_2)   \\[.6ex]
          \widetilde{\beta}_1^{\,0}\ne 0,\,
             \widetilde{\beta}_2^{\,0}\ne 0 \\[.6ex]
          p_{1\ast}\widetilde{\beta}_1^{\,0}
            + p_{2\ast}\widetilde{\beta}_2^{\,0}=\beta
        \end{array}
                }
         \left\{
          \begin{array}{l}
            \eta=(\Gamma_1,\Gamma_2,I) \\
             \mbox{admissible} \\
             \mbox{triple}     \\
             \mbox{for $Y_1\cup_ E Y_2$}
          \end{array}
          \left|
            \begin{array}{l}             
              \bullet\hspace{1ex}
                  b(\Gamma_1)=\widetilde{\beta}_1^{\,0}+l_1\gamma,\,
                  b(\Gamma_2)=\widetilde{\beta}_2^{\,0}+l_2\gamma\,, \\
              \hspace{2ex} l_1+l_2 =E\cdot \widetilde{\beta}_1^{\,0}\,,\;
                           l_1,\,l_2\in{\footnotesizeBbb Z}_{\ge 0}\,;
                                                                 \\[.6ex]
              \bullet\hspace{1ex}
               g(\eta) = g\,,\;
               k_1+k_2 = k\,; \\[.6ex]
              \bullet\hspace{1ex}
               \sum_{i}\mu_{1,i}\,=\,l_2\,;\\[.6ex]
              \bullet\hspace{1ex}
               I\subset \{1,\,\ldots,\, k\}\,,\; |I|=k_1\,.
            \end{array}
          \right.
         \right\}\;
 $$
 $$
  \coprod
  \coprod_{\tiny
           \begin{array}{c}
             0\ne\widetilde{\beta}_1^{\,0}\in NE(Y_1) \\[.6ex]
             E\cdot \widetilde{\beta}_1^{\,0}\ge 0     \\[.6ex]
             p_{1\ast}\widetilde{\beta}_1^{\,0}=\beta
           \end{array}
          }
   \left\{
          \begin{array}{l}
            \Gamma_1 :\; \mbox{admissible}   \\
             \mbox{weighted graph}           \\
             \mbox{for $(Y_1,E)$}
          \end{array}
      \left|
            \begin{array}{l}
              \bullet\hspace{1ex}
               b(\Gamma_1)
                =\widetilde{\beta}_1^{\,0}
                  +(E\cdot\widetilde{\beta}_1^{\,0})\gamma\,;\\[.6ex]
              \bullet\hspace{1ex}
               g(\Gamma_1) = g\,,\;
               \mbox{$k$-many legs}\,; \\[.6ex]
              \bullet\hspace{1ex}
               \mbox{no roots}\,.      \\[.6ex]
            \end{array}
      \right.
   \right\}
 $$
 $$
  \hspace{2ex}
  \coprod
  \coprod_{\tiny
           \begin{array}{c}
             0\ne \widetilde{\beta}_2^{\,0}\in NE(Y_2) \\[.6ex]
             p_{2\ast}\widetilde{\beta}_2^{\,0}=\beta
           \end{array}
          }
   \left\{
          \begin{array}{l}
            \Gamma_2 :\; \mbox{admissible}   \\
             \mbox{weighted graph}           \\
             \mbox{for $(Y_2,E)$}
          \end{array}
      \left|
            \begin{array}{l}
              \bullet\hspace{1ex}
               b(\Gamma_2)
                =\widetilde{\beta}_2^{\,0}\,;\\[.6ex]
              \bullet\hspace{1ex}
               g(\Gamma_2) = g\,,\;
               \mbox{$k$-many legs}\,; \\[.6ex]
              \bullet\hspace{1ex}
               \mbox{no roots}\,.      \\[.6ex]
            \end{array}
      \right.
   \right\}\,.
 $$
} % end-footnotesize

\noindent
(Note that in the above expression, up to the orderings of legs,
 the whole set in respectively line 2 and line 3 of the equation
 is either the empty set or a singleton.)

We now summarize the discussions in Sec.\ 2 and Sec.\ 3
 in the following theorem:

\bigskip

\noindent
{\bf Theorem 3.3 [J.\ Li's degeneration formula refined for blow-up].}
{\it
 Let $Z$ be a smooth subvariety of a smooth variety $X$. Consider
  the family $W\rightarrow {\Bbb A}^1$ by blowing $X\times{\Bbb A}^1$
  along $Z\times{\mathbf 0}$.
 Let
  $\beta\in A_1(X)$ and
  $\Psi_{(g,k;\beta)}^X(\,\cdot\,,\,\cdot\,)$
   be the usual Gromov-Witten invariants of $X$ defined via the stack
    $\overline{\cal M}_{g,k}(X,\beta)$ of stable maps of genus $g$,
     $k$ marked points, and to the curve class $\beta$ in $X$.
 Then, with the same notation as in Corollary 2.2,
   %\Corollary ??? [J.\ Li's degeneration formula]
  there is a canonical choice $\Omega_{(g,k;\beta)}$ of a finite set
  of admissible triples that depends only on ${\cal N}_{Z/X}$,
  $(g,k)$, and $\beta\in A_1(X)$
  {\rm (}more precisely, the semi-group $\NE(X)_{\scriptsizeBbb Z}${\rm )}
  such that
 $$
  \Psi_{(g,k;\beta)}^X(\alpha,\zeta)\;
   =\; \sum_{\eta\in\overline{\Omega}_{(g,k;\beta)}}\;
        \frac{{\mathbf m}(\eta)}{|\Eq(\eta)|}\,
          \sum_{j\in K_{\eta}}\,
           \left[ \Psi_{\Gamma_1}^{Y_1^{\rel}}(j_1^{\ast}\alpha(0),
                    \zeta_{\eta,1,j})\,
                  \bullet\,
                  \Psi_{\Gamma_2}^{Y_2^{\rel}}(j_2^{\ast}\alpha(0),
                    \zeta_{\eta,2,j})
           \right]_0\,.
 $$
Similarly for the associated cycle form of the formula.
The same statements hold with $\beta\in A_1(X)$ replaced by 
 $\beta\in H_2(X;{\Bbb Z})$.
} % end-corollary

\bigskip

\noindent
{\it Remark 3.4 {\rm [}generalization{\rm ]}.}
The same technique and a similar statement apply also for
 a degeneration $W\rightarrow {\Bbb A}^1$ from a sequence of blow-ups
 of a trivial family $X\times {\Bbb A}^1$ and for the situation when
 the smooth subscheme $Z\subset X$ to be blown up is disconnected,
 with each connected component of its own dimension.
For more general degenerations, one needs to understand better the
 monodromy effect around a singular fiber of $W/{\Bbb A}^1$ and
 how $\NE(W/{\Bbb A}^1)_{\scriptsizeBbb Z}$ behaves over ${\Bbb A}^1$.

\vspace{2em}
%references
%endfootnotesize

\end{document}